\documentclass[11pt]{amsart}

\usepackage{amsfonts,amsmath,latexsym,amssymb,verbatim,amsbsy,amsthm}
\usepackage{graphicx}
\usepackage{caption}
\usepackage{float}
\usepackage{wrapfig}
\usepackage{mathrsfs}

\usepackage[top=1in, bottom=1in, left=1in, right=1in]{geometry}

\usepackage[dvipsnames]{xcolor}

\usepackage[colorlinks=true, pdfstartview=FitV, linkcolor=RoyalBlue,citecolor=ForestGreen, urlcolor=blue]{hyperref}

\theoremstyle{plain}
\newtheorem{THEOREM}{Theorem}[section]

\newtheorem{theorem}[THEOREM]{Theorem}

\theoremstyle{definition}

\theoremstyle{remark}

\DeclareMathOperator{\Supp}{Supp} %
\DeclareMathOperator{\supp}{supp} %
\DeclareMathOperator{\Tr}{Tr} %
\def \a {\alpha}

\def \d {\delta}

\def \e {\varepsilon}

\def \n {\nabla}

\def \th {\theta}

\def \D {\Delta}

\def \L {\Lambda}

\def \bu {{\bf u}}

\def \bx {{\bf x}}

\def \cA {\mathcal{A}}

\def \cC {\mathcal{C}}
\def \cD {\mathcal{D}}

\def \cI {\mathcal{I}}

\def \cM {\mathcal{M}}

\def \cP {\mathcal{P}}

\def \dH {\dot{H}}\def \dW {\dot{W}}
\def \rmin{\underline{\rho}}

\newcommand{\Z}{\ensuremath{\mathbb{Z}}}   
\newcommand{\R}{\ensuremath{\mathbb{R}}}   
\newcommand{\T}{\ensuremath{\mathbb{T}}}   

\def \lan {\langle}
\def \ran {\rangle}

\def \p {\partial}

\def \ss {\subset}

\def \dx  {\, \mbox{d}x}

\def \dy  {\, \mbox{d}y}
\def \dz  {\, \mbox{d}z}
\def \dr  {\, \mbox{d}r}
\def \ds  {\, \mbox{d}s}
\def \dw  {\, \mbox{d}w}

\def \ddt  {\frac{\mbox{d\,\,}}{\mbox{d}t}}
\def \DDt  {\frac{\mbox{D\,\,}}{\mbox{D}t}}

\begin{document}

\title[Unidirectional flocks: Singular models]{ Unidirectional flocks in hydrodynamic Euler Alignment system II:  singular models.}

\author{Daniel Lear} \author{Roman Shvydkoy}

\address{Department of Mathematics, Statistics, and Computer Science, University of Illinois at Chicago}

\email{lear@uic.edu}
\email{shvydkoy@uic.edu}

\date{\today}



\thanks{\textbf{Acknowledgment.} Research of RS is supported in part by NSF
	grant DMS-1813351}

\begin{abstract}
In this note we continue our study of unidirectional solutions to  hydrodynamic Euler alignment systems
with strongly singular communication kernels $\phi(x):=|x|^{-(n+\a)}$ for $\a\in(0,2)$. The solutions describe unidirectional parallel motion of agents governing multi-dimensional collective behavior of flocks. Here, we consider the range $1<\a<2$ and establish the global regularity of smooth solutions, together with a full description of their long time dynamics. Specifically, we develop the flocking theory of these solutions and show long time convergence to traveling wave with rapidly aligned velocity field. 
\end{abstract}

\maketitle

\section{Introduction and statement of main results}
We consider the following  hydrodynamic Euler Alignment System for density $\rho(x,t)$ and velocity $\bu(x,t)=(u^{1}(x,t),\ldots,u^{n}(x,t))$  :
\begin{equation}\label{e:CSHydro}
(x,t)\in\mathbb{R}^{n}\times\mathbb{R}^{+}\qquad \left\{
\begin{split}
\partial_t \rho +\nabla\cdot (\rho \bu )&= 0, \\
\partial_t  \bu +\bu \cdot\nabla \bu &= \mathcal{L}_{\phi} (\rho  \bu ) -  \mathcal{L}_{\phi}(\rho)\bu,
\end{split}\right.
\end{equation}
subject to initial condition
$$\left( \rho(\cdot,t), \bu (\cdot,t)\right)|_{t=0}=(\rho_{0},\bu_{0}).$$

The system (\ref{e:CSHydro}) arises as a macroscopic realization of the Cucker-Smale agent-based dynamics \cite{CS2007a,CS2007b}, which describes collective motion of $N$ agents adjusting their directions to a weighted average of velocities of its neighbors:
\begin{equation*}
(\bx_{i},\mathbf{v}_{i})\in\mathbb{R}^{n}\times\mathbb{R}^{n}\qquad \left\{
\begin{array}{@{} l l @{}}
\dot{\bx}_{i}&\hspace{-0.2 cm}= \mathbf{v}_{i}, \\
\dot{\mathbf{v}}_{i}&\hspace{-0.2 cm}=\frac{1}{N}\sum_{j=1}^{N}\phi(|\bx_{i}-\bx_{j}|)(\mathbf{v}_{j}-\mathbf{v}_{i}).
\end{array}\right. 
\end{equation*}

We refer to \cite{HL2009,FK2019,HT2008} for full details and  rigorous derivations. Typical assumptions on $\phi(r)$ include monotonic decay at infinity and non-degeneracy, $\phi(r) >0$,  thus reflecting the intuition that alignment becomes weaker, yet persistent, as the distance becomes larger.  When communication remains sufficiently strong at infinity, expressed by the ``fat tail" condition
\begin{equation}\label{e:ft}
	\int_0^\infty \phi(r) \dr = \infty,
\end{equation}
the system \eqref{e:CSHydro} (as well as its discrete counterpart) exhibits alignment dynamics, that is for any global strong solution,  
\begin{equation*}
\mathcal{A}(t):=\underset{\left\lbrace x, y \right\rbrace\in\Supp \rho(\cdot,t)}{\text{max}}|\bu(x,t)-\bu( y ,t)|\to 0\qquad \text{as} \quad t\to\infty
\end{equation*}
exponentially fast, and the diameter of the flock remains globally bounded:
$$
\mathcal{D}(t)\leq \bar{\cD}<\infty \qquad \text{where }\quad \mathcal{D}(t):=\underset{\left\lbrace x, y \right\rbrace\in\Supp \rho(\cdot,t)}{\text{max}}|x- y |.
$$

For models with singular kernels given by $\phi(x)=|x|^{-(n+\a)}$  for  $0<\alpha<2$ the
operator $\mathcal{L}_{\phi}\equiv\mathcal{L}_{\a}$ becomes the (negative of) classical fractional Laplacian:
\[
\mathcal{L}_{\a}(f)(x)=-\Lambda_{\a}(f)(x)\equiv \text{p.v.}\int_{\R^n}\frac{f(y)-f(x)}{|x-y|^{n+\a}}\, \dy    \qquad \Lambda_{\a}:=(-\Delta)^{\a/2},   \qquad 0<\a<2.
\]
The corresponding  alignment term on the right hand side of the momentum equation in \eqref{e:CSHydro} is then given by the following singular integral:
\begin{equation}\label{e:alignmentalpha} 
\mathcal{C}_{\a}(\bu,\rho)= -\Lambda_{\a} (\rho  \bu ) +  \Lambda_{\a}(\rho)\bu=\text{p.v.}\int_{\R^n}\frac{\bu(y)-\bu(x)}{|x-y|^{n+\a}} \rho(y)\, \dy.   
\end{equation}
In view of no-vacuum condition $(\rho_0>0)$ necessary to develop a well-posedness theory we consider the periodic domain $\T^n$, where a uniform lower bound on the density is compatible with finite mass.
When dealing with the $n$-dimensional torus, the term \eqref{e:alignmentalpha} can be expressed in terms of the periodized kernel
\[
\phi_{\a}(z):=\sum_{k\in \Z^n}\frac{1}{|z+2\pi k|^{n+\a}}, \qquad 0<\a<2,
\]
which preserve the essential long range but less dominant interactions. In the rest of the paper, we assume that $\bu(\cdot,t)|_{\T^n}$ and likewise $\rho(\cdot,t)|_{\T^n}$ are extended periodically onto the whole space $\R^n$. The alignment term \eqref{e:alignmentalpha} then becomes a fractional elliptic operator:
\[
\mathcal{C}_{\a}(\bu,\rho)=\text{p.v.}\int_{\R^n}\left(\bu(x+z)-\bu(x)\right) \rho(x+z)\frac{\dz}{|z|^{n+\a}}=\text{p.v.}\int_{\T^n}\left(\bu(x+z)-\bu(x)\right) \rho(x+z) \phi_{\a}(z) \dz,   
\]
with the density controlling uniform ellipticity. Written in this form,  system \eqref{e:CSHydro} resembles the fractional Burgers equation with non-local non-homogeneous dissipation.\\ 

In \cite{ST1,ST2,ST3} Tadmor and the second author proved global existence of smooth solutions for the one-dimensional  system \eqref{e:CSHydro} with alignment term given by  \eqref{e:alignmentalpha} in the  full range  $0< \a < 2$, with focus on the most difficult critical case $\a=1$.   In addition, the authors proved in \cite{ST2} that all regular solutions converge exponentially fast to a so called \emph{flocking state}, consisting of a traveling wave, $\bar{\rho}(x,t)=\rho_\infty(x-t\bar{u})$, with a fixed speed $\bar{u}$,
\begin{equation*}
| u(\cdot,t)- \bar{u}|_X + |\rho(\cdot,t) - \bar{\rho}(\cdot,t)|_Y  \xrightarrow{t\rightarrow\infty}0, \qquad  \bar{u}:=\frac{\cP_0}{\cM_0}.
\end{equation*}
Here the average velocity, $\overline{u}$, is dictated by the conserved mass and momentum, 
\[
\cM_0=\int_\T \rho_0(x) \dx, \quad \cP_0=\int_\T (\rho_0u_0)(x) \dx.
\] 
Parallel to these works, Do et.al. in \cite{DKRT2018} treated the case $0<\a<1$, where they proved global existence with the use of the modulus of continuity method as in Kiselev et. al. \cite{KNV2007}. In either approach the problem requires utilization of refined tools from regularity theory of fractional parabolic equations, and reduces to verification of  a continuation criterion either in terms of  $u_x \in L^1\left([0,T_0);L^{\infty}\right)$, \cite{ST1}, or in terms of $\rho_x \in L^1\left([0,T_0);L^{\infty}\right)$, \cite{DKRT2018}.

Global well-posedness theory for these singular models has been developed only in 1D mainly due to presence of an additional conserved quantity
\begin{equation}\label{e:entropy1D}
e:=u_x -\Lambda_{\a}(\rho), \qquad e_t +(eu)_x=0,
\end{equation}
Thanks to these relations, one can compare the regularity of $u$ and $\rho$ and, using the compactness of the 1D torus, obtain a global-in-time positive lower bound on the density, thanks to which the ``good'' term on the right-hand side of \eqref{e:alignmentalpha} does not disappear. That method unfortunately fails in higher dimension. 
In multi-dimensional settings  the corresponding quantity is given by
\[
e:= \nabla \cdot \bu -\Lambda_{\a}(\rho)
\]
and satisfies
\[
e_t+\nabla \cdot(\bu e)=\left(\nabla \cdot \bu \right)^2 -\text{Tr}[\left(\nabla \bu\right)^2].
\]
Lack of control on $e$ in this case is part of the reason why in multiple dimensions the model has no developed regularity theory. The two exceptions are small initial data results proved in \cite{Shv2018} and \nolinebreak \cite{DMPW2019}.

Although in 1D local existence for singular models appeared in the first papers by Shvydkoy, Tadmor \cite{ST1} for $\a\geq 1$, and Do et. al. \cite{DKRT2018} for $0 < \a < 1$, it was not properly addressed in higher dimensions, with a proper continuation criterion. We fill this gap with the following result.
\begin{theorem}[Local existence of classical solutions]\label{t:lwpsing}
	Suppose  $m > \frac{n}{2}+1$, $0<\a<2$, and 
	\[
	(\bu_0,\rho_0) \in H^{m+1}(\T^n) \times H^{m+\a}(\T^n),
	\]
	with $\rho_0(x) > 0$ for all $x\in \T^n$. Then there exists time $T_0>0$  and a unique non-vacuous solution to \eqref{e:CSHydro} on time interval $[0,T_0)$ in the class
	\begin{equation*}
	\bu \in C_w([0,T_0); H^{m+1}) \cap  L^2([0,T_0); \dH^{m+1+ \a/2} ), \quad \rho \in C_w([0,T_0); H^{m+\a}) .
	\end{equation*}
	Moreover, any such solution satisfying 
	\begin{equation}\label{e:BKMsing}
	\sup_{t \in [0,T_0)} (|\n \rho(t)|_\infty+ | \n \bu(t) |_\infty )<\infty
	\end{equation} 
	can be extended beyond $T_0$.
\end{theorem}

\begin{wrapfigure}{l}{0.38\textwidth}
	\centering
     \includegraphics[width=0.35\textwidth]{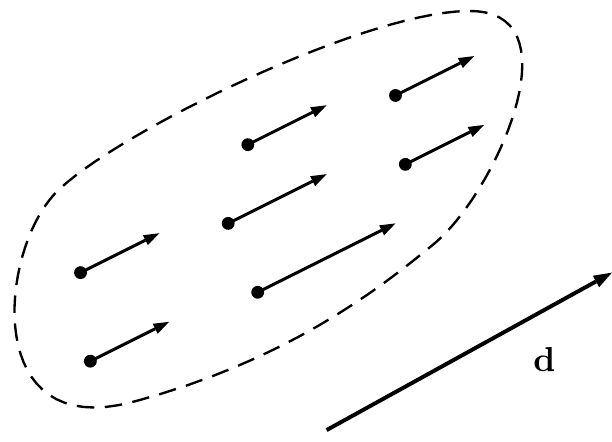}
     \vspace*{0.5 cm}
\end{wrapfigure}

\noindent
One class of solutions that behaves like 1D is the class of unidirectional oriented flows introduced in \cite{LS-uni1} for the case of smooth kernels. These are given by
\begin{equation*}\label{e:ansatz}
\hspace{6 cm}\bu(x,t)=u(x,t)\mathbf{d}, \quad \mathbf{d}\in \mathbb{S}^{n-1}, \quad u:\R^{n}\times \R^{+}\rightarrow \R.
\end{equation*}
The same conservation law \eqref{e:entropy1D} holds for the entropy
\[
\hspace{6.5 cm} e:=\mathbf{d}\cdot \nabla u -\Lambda_{\a}(\rho), \qquad \p_t e + \mathbf{d}\cdot\nabla (u e)=0,
\]
although in this case the entropy does not control the full gradient of the velocity.  Let us make a couple of remarks about the unidirectional ansatz itself. 
First, by the maximum principle of the velocity field applied in any direction perpendicular to $\mathbf{d}$ one can see that the ansatz is preserved in time. Second, in view of rotational invariance of the Euler Alignment System, we can postulate that \textbf{d} points in the direction of the $x_1$-axis. So, we can assume
\[
\bu(x,t):=\lan u(x,t),0,\ldots,0 \ran \qquad \text{for} \quad u:\R^{n}\times \R^{+}\rightarrow \R.
\]
Note  that  the  non-trivial  component $u(x,t)$  may  depend  on  all  coordinates.   So,  our  solutions exhibits features of a 1D flow, yet being on $\R^n$ represent solutions of a multi-D system of scalar conservation laws:
\begin{equation}\label{e:CSHansatz}
(x,t)\in\mathbb{R}^n\times\mathbb{R}^{+}\qquad \left\{
\begin{array}{@{} l l @{}}
\partial_t \rho  +\partial_1(\rho u)&\hspace{-0.3 cm}=0, \\
\partial_t u+ \frac12\partial_1 (u^2)&\hspace{-0.3 cm}=\cC_{\a}(u,\rho).
\end{array}\right. 
\end{equation}
In this paper we will continue the study of unidirectional solutions initiated in \cite{LS-uni1}, now in the context of singular models. As in \cite{ST1,ST2,ST3}, our methodology will be to extract \emph{quantitative enhancement} estimates for the dissipation term, using an adaptation of the non-linear maximum principle as in Constantin and Vicol's proof for the critical SQG \cite{CV2012}, that  yields global existence and, moreover, allows us to  completely describe the long time behavior --- exponential convergence towards a flocking state. The main result  summarized in the following theorem covers the global regularity and flocking behavior for singular kernels in the range $1<\a<2$.

\begin{theorem}\label{t:flock-singular}
Suppose $m \geq 3$ and $1<\a<2$. Let $(u_0,\rho_0) \in H^{m+1}(\T^n) \times H^{m+\a}(\T^n)$ with $\rho_0(x) > 0$ for all $x\in \T^n$. Then there exists a unique non-vacuous global in time solution to \eqref{e:CSHansatz} in the class
\[
u \in C_w([0,\infty); H^{m+1}) \cap  L^2([0,\infty); \dH^{m+1+ \a/2} ), \quad \rho \in C_w([0,\infty); H^{m+\a}).
\]
Moreover, the solution obeys uniform bounds on the density
\begin{equation}\label{e:rhopositive}
c_0 \leq \rho(x,t)\leq C_0, \quad t\geq 0,
\end{equation}
and strong flocking: $\exists \bar{\rho}\in H^{m+\a}$ such that
\begin{equation*}
\|u(t)-\bar{u}\|_{W^{2,\infty}} +\|\rho(\cdot, t) - \bar{\rho}(\cdot -\bar{u}t) \|_{C^{\gamma}} \leq C e^{-\d t}, \qquad t>0,\quad  (0<\gamma<1). 
\end{equation*}
\end{theorem}
As before, the limiting velocity $\bar{u}$ is determined from the initial conditions due to conservation of mass and  momentum.

There does not seem to be a rule in either 1D or our situation  on how to determine the limiting density distribution of the flock  $\bar{\rho}$ -- this appears to be an \emph{emerging} quantity of the dynamics. However, the entropy estimates done in \cite{LS-entropy} show that, at least on the periodic domain the size of $e$ directly controls how far $\bar{\rho}$ is from the uniform distribution. \\

\noindent
\textit{Notation:}
For convenience, to avoid clutter in computations, function arguments (time and space) will be omitted whenever they are obvious from context. Moreover, we use the notation $f\lesssim g$ when there exists a constant $C >0$ independent of the parameters of interest such that $f\leq C g$. We also use $|\cdot|_p$, $1\leq p\leq \infty$, to denote the classical $L^p$-norms, and $\| \cdot \|_X$ to denote all other norms.\\

\noindent
\textit{Organization:}
In Section \ref{s:LWP} we prove a local existence and a continuation criterion result in Sobolev spaces with minimal requirements needed for what follows.  In Section \ref{s:GWP}, as a direct application of the continuation criterion, we obtain a global existence result for unidirectional parallel motion, and provide higher order control estimates on solutions to prove a strong flocking result.

\section{Local well-posedness and continuation criteria}\label{s:LWP}

We will be casting our regularity theory for singular models on the periodic domain $\T^n$ and for non-vacuous solutions only. This is motivated by technical reasons rather than applications, although one can argue that periodic  conditions are suitable for studying  flocks in the bulk. The primary reason is that we require uniform parabolicity of the commutator \eqref{e:alignmentalpha} for estimates to go through. Such parabolicity depends on the pointwise bound $\rho >c_0>0$, which is consistent with finite mass of the flock only on bounded domains.  

Necessity of the no-vacuum condition can be easily seen by the following example in 1D. Let us consider
a local kernel for simplicity, $\supp \phi \ss B_{1}(0)$. Let initial density $\supp \rho_0 \ss B_{\e}(0)$, while $u_0 = 1$ on $B_{10}(0)$, $u_0 = 0$ on $\R\setminus B_{10+\e}(0)$ and smooth in between. Then the density will remain in $B_2(0)$ for a time period of at least $t<1$, due to $u\leq 1$. During this time the momentum equation will remain pure Burgers, hence the solution will evolve into a shock at a time $t \sim \e<1$.  
A more subtle blowup can be constructed even  for a global singular kernel on the periodic domain  when the density vanishes at just one point \cite{AC2019}. Earlier \cite{Tan2017} demonstrated growth of $\|\rho\|_{C^1}$ as $t\to \infty $ for a similar density configuration.\\

Performing energy estimates in the same fashion as for smooth models \cite{LS-uni1} will inevitably create a derivative overload on the density. Instead we consider another ``almost conserved'' quantity 
\begin{equation*}
e= \n \cdot \bu -\Lambda_{\a}(\rho),
\end{equation*}
which satisfies the equation
\begin{equation}\label{e:eq}
e_t + \nabla \cdot(\bu e)  = (\n \cdot \bu)^2 -  \Tr(\n \bu)^2.
\end{equation}
Let us derive it in general for the sake of completeness. 
Since $\phi_{\a}$ is a convolution kernel, we have 
\begin{equation}\label{aux_1}
\partial_t \Lambda_{\a} (\rho) +\n \cdot \Lambda_{\a} (\rho \bu)=0.
\end{equation}
Taking the divergence of the velocity equation, we obtain
\begin{equation}\label{aux_2}
\partial_t (\n \cdot \bu)+\n \cdot [ \bu \cdot\nabla \bu ]=\n \cdot [ \bu\, \Lambda_{\a}(\rho) ]-\n \cdot \Lambda_{\a}(\rho  \bu ) 
\end{equation}
with
$$\n \cdot [ \bu\, \Lambda_{\a} (\rho) ]=\Lambda_{\a} (\rho) \n \cdot \bu+ \bu \cdot \nabla \Lambda_{\a} (\rho)$$
and
$$\n \cdot [ \bu \cdot\nabla  \bu ]=\Tr (\nabla  \bu )^{2} + \bu \cdot \nabla(\n \cdot \bu) .$$
On one hand, combining (\ref{aux_1}) and (\ref{aux_2}), we obtain that
\begin{equation*}
\partial_t e -\Lambda_{\a} (\rho) \n \cdot \bu + \bu \cdot\nabla e +\Tr (\nabla  \bu )^{2} =0.
\end{equation*}
Adding and subtracting now  $(\n \cdot \bu )^{2}$ produces \eqref{e:eq}. It is clear that in 1D the right hand side of \eqref{e:eq} vanishes, and one obtains a perfect continuity law.

From the order of terms that enter into the formula for $e$, it is clear that the natural correspondence in regularity for state variables involved is  $(\bu\in H^{m+1}) \sim (\rho \in H^{m+\a})$. The grand quantity to be estimated is
\[
Y_m = \|\bu\|_{H^{m+1}}^2 + \|e\|_{H^m}^2 + |e|_\infty +  |\rho|_1 + |\rho^{-1}|_\infty,
\]
which is equivalent to $Y_m \sim  \|\bu\|_{H^{m+1}}^2 + \|\rho\|_{H^{m+\a}}^2+ |\rho^{-1}|_\infty$ thanks to
the  coercivity estimate for $\Lambda_{\a}$:
\begin{equation}\label{e:coer}
c_1 \|f\|_{\dot{H}^{\a}}-c_2 |f|_{2}\leq |\Lambda_{\a}f|_2 \leq c_3 \|f\|_{\dot{H}^{\a}}-c_4 |f|_{2}
\end{equation}

Our strategy will be very similar to the smooth case \cite{LS-uni1}, where we obtain local solutions via viscous regularization, and prove a continuation criterion via a priori estimates on $Y_m$. 

To actually produce local solutions we consider viscous regularization of the system \eqref{e:CSHydro} and we assume throughout that $m > \frac{n}{2} + 1$ and $0<\a<2$.
\begin{equation}\label{e:CSHvisc}
\left\{
\begin{split}
\p_t\rho + \n \cdot (\rho \bu) & =  \e \D \rho, \\
\p_t\bu + \bu \cdot \n \bu &= \cC_{\a}(\bu,\rho)+ \e \D \bu.
\end{split}\right.
\end{equation}
So, let us start with \eqref{e:CSHvisc} and consider the mild formulation
\begin{equation*}
\begin{split}
\rho(t) &= e^{\e t \D} \rho_0 - \int_0^t e^{\e(t-s) \D} \n \cdot (\bu \rho)(s) \ds \\
\bu(t) &= e^{\e t \D} \bu_0 - \int_0^t e^{\e(t-s) \D} \bu \cdot \n \bu (s) \ds + \int_0^t e^{\e(t-s) \D}  \cC_{\a}(\bu,\rho)(s) \ds.
\end{split}
\end{equation*}
Let us denote by $Z = (\rho,\bu)$ the state variable of our system and by $T[Z](t)$ the right hand side of the mild formulation. 
In order to apply the standard fixed point argument we have to show that $T$ leaves the set $C([0,T_{\d,\e}); B_\d(Z_0))$ invariant, where $B_\d(Z_0)$ is the ball of radius $\e$ around initial condition $Z_0$, and that it is a contraction. We limit ourselves to showing details for invariance as the estimates involved in proving Lipschitzness are similar. 

First we assume that $\rho$ has no vacuum: $\rho_0(x) \geq  c_0>0$. Since the metric we are using for $\rho\in H^{m+\a}$ controls $L^\infty$ norm, if $\d>0$ is small enough then for any $\| \rho - \rho_0\|_{H^{m+\a}} < \d$ one obtains 
\begin{equation*}
\rho(x) > \frac12 c_0.
\end{equation*}
So, let us assume that $Z \in C([0,T); B_\d(Z_0))$. It is clear that $\|e^{\e t\D} Z_0 - Z_0\|< \frac{\d}{2}$ provided time $t$ is short enough. The $Z$ has some bound $\|Z\| \leq C$.  Using that let us estimate the norms under the integrals. First, recall that $ \| \L_\a e^{\e t\D}\|_{L^2 \to L^2} \lesssim \frac{1}{(\e t)^{\a/2}}$. In the case $\a \geq 1$, we have
\begin{multline*}
\left| \p^m \L_\a \int_0^t e^{\e(t-s) \D} \n \cdot (\bu \rho)(s) \ds \right|_2 \lesssim \int_0^t \frac{1}{(t-s)^{\a/2}} | \p^{m+1} (\bu\rho)(s) |_2 \ds \\
\leq \int_0^t \frac{1}{(t-s)^{\a/2}}  \|\bu\|_{\dH^{m+1}} \|\rho\|_{\dH^{m+\a}} \ds \leq C^2 t^{1-\a/2} < \frac{\d}{2},
\end{multline*}
provided $T = T(\d,\e)$ is small enough. In the case $\a<1$, we combine instead one full derivatives with the heat semigroup, and the rest $\p^{m+\a}$ gets applied to $\bu\rho$, which produces a similar bound. 

Moving on to the $\bu$-equation, we have
\begin{multline*}
\left| \p^{m+1} \int_0^t e^{\e(t-s) \D}  \bu \cdot \n \bu (s)  \ds \right|_2 \lesssim \int_0^t \frac{1}{(t-s)^{1/2}} | \p^{m} (\bu \cdot \n \bu)(s) |_2 \ds \\
\leq \int_0^t \frac{1}{(t-s)^{\a/2}}  \|\bu\|_{\dH^{m+1}} \|\bu\|_{\dH^{m}}\ds \leq C^2 t^{1/2} < \frac{\d}{4}.
\end{multline*}
As to the commutator form, for $\a\leq 1$ the computation is very similar: we combine one derivative with the heat semigroup and for the rest we use \eqref{e:coer}:
\[
| \p^m \cC_{\a}(\bu,\rho)|_2 \lesssim  \|\bu\|_{m+\a} \|\rho\|_{m+\a}  < C^{2},
\]
and the rest follows as before.  When $\a > 1$ we combine $\a$ derivatives with the semigroup, and the rest follows as before. 

We have proved that $\|T[Z](t) - Z_0 \| <\d$, for a short time and hence, $T$ leaves $C([0,T(\d,\e)); B_\d(Z_0))$ invariant. 
The obtained interval of existence of course depends on $\e$ as it enters into all the estimates of the integrals.  In order to conclude the local existence argument we still have to show that our apriori bound is independent of $\e$.  This would allow us to extend $T_{\e,\delta}$ to a time dependent on the initial condition only. Then the classical compactness argument would apply to pass to the limit as $\e\to 0$ in the same state space $C ([0,T); H^{m+\a}\times H^{m+1}).$\\

Now let us make a priori estimates for viscous solutions independent of $\e$. Note that the dissipation terms in all the following computations are negative and as such will be ignored.

First, evaluating the continuity equation at a point of minimum $x_-$ and denoting $\rmin = \min \rho$ we readily obtain
\[
\ddt \rmin = - \rmin \n \bu + \e \D \rho (x_-) \geq - \rmin |\n \bu|_\infty.
\] 
Hence,
\[
\ddt | \rho^{-1}|_\infty \leq | \rho^{-1}|_\infty |\n \bu|_\infty \leq  |\n \bu|_\infty Y_m.
\]
Furthermore,
\begin{equation}\label{e:euu}
\ddt |e|_\infty \leq |\n \bu|_\infty |e|_\infty + |\n \bu|_\infty^2 \leq |\n \bu|_\infty Y_m.
\end{equation}
Let us continue  with estimates on the $e$-quantity. We have (dropping integral signs)
\[
\ddt \|e \|_{\dH^m}^2 \leq  \p^m e \bu \cdot \n \p^m e + \p^m e [ \p^m( \bu \cdot \n e) - \bu \cdot \n \p^m e] + \p^m e  \p^m(e \n \cdot \bu) + \p^m e [(\n \cdot \bu)^2 -  \Tr(\n \bu)^2]
\] 
In the first term we integrate by part and estimate
\[
| \p^m e \bu \cdot \n \p^m e| \leq \|e \|_{\dH^m}^2 |\n \bu|_\infty.
\]
For the next term we use the classical commutator estimate
\begin{equation}\label{e:classcomm}
|\p^k(fg)-f\p^k g|_2 \leq |\nabla f|_{\infty}\|g\|_{\dot{H}^{k-1}}+\|f\|_{\dot{H}^k}|g|_{\infty}
\end{equation}
to obtain that
\[
|\p^m e [ \p^m( \bu \cdot \n e) - \bu \cdot \n \p^m e]| \leq \|e \|^2_{\dH^m} |\n \bu|_\infty +\|e \|_{\dH^m} \|\bu  \|_{\dH^m} |\n e|_\infty.
\]
Using Gagliardo-Nirenberg inequality we estimate the latter term as
\[
\|e \|_{\dH^m} \|\bu  \|_{\dH^m} |\n e|_\infty \leq \|e \|_{\dH^m} \|\bu\|_{\dH^{m+1}}^{\th_1} |\n \bu|_\infty^{1-\th_1} \| e\|_{\dH^m}^{\th_2} |e|_\infty^{1-\th_2},
\]
where $\th_1 = \frac{n-2(m-1)}{n-2m}$ and $\th_2 = \frac{2}{2m-n}$. The two exponents add up to $1$, so by the generalized Young inequality,
\[
\leq (\| e\|_{\dH^m}^2 + \|\bu\|_{\dH^{m+1}}^2) ( |e|_\infty +  |\n \bu|_\infty) \leq ( |e|_\infty +  |\n \bu|_\infty) Y_m
\]
Next term in the $e$-equation is estimated by the product formula
\begin{equation*}
| \p^m (f g) |_2 \leq \| f\|_{H^m} |g|_\infty + |f|_\infty \|g\|_{H^m}.
\end{equation*}
So, we have
\[
| \p^m e  \p^m(e \n \cdot \bu) | \leq \|e \|_{\dH^m}^2 |\n \bu|_\infty + \|e \|_{\dH^m} |e|_\infty  \|\bu\|_{\dH^{m+1}} \leq ( |e|_\infty +  |\n \bu|_\infty) Y_m.
\]
Finally,
\[
| \p^m e [(\n \cdot \bu)^2 -  \Tr(\n \bu)^2] | \leq \|e \|_{\dH^m} \|\bu\|_{\dH^{m+1}}  |\n \bu|_\infty \leq  |\n \bu|_\infty Y_m.
\]
Thus, 
\begin{equation*}
\ddt \|e \|_{\dH^m}^2 \leq ( |e|_\infty +  |\n \bu|_\infty) Y_m.
\end{equation*}

Next perform the main technical estimate on the velocity equation. We have
\[
\p_t \| \bu\|_{\dH^{m+1}}^2  = - \p^{m+1} (\bu \cdot \n \bu) \cdot \p^{m+1}  \bu +  \p^{m+1} \cC_{\a}(\bu,\rho)\cdot \p^{m+1}  \bu.
\]
The transport term is estimated using the classical commutator estimate
\[
\p^{m+1} (\bu \cdot \n \bu) \cdot \p^{m+1}  \bu = \bu \cdot \n  (\p^{m+1} \bu) \cdot \p^{m+1}  \bu + [\p^{m+1} , \bu] \n \bu \cdot \p^{m+1}  \bu 
\]
Then
\[
\bu \cdot \n  (\p^{m+1} \bu) \cdot \p^{m+1}  \bu = - \frac12 (\n \cdot \bu) | \p^{m+1}  \bu|^2 \leq |\n \bu|_\infty \| \bu\|_{\dH^{m+1}}^2,
\]
and using  \eqref{e:classcomm}  we obtain
\[
|[\p^{m+1} , \bu] \n \bu \cdot \p^{m+1}  \bu  | \leq  |\n \bu|_\infty \| \bu\|_{\dH^{m+1}}^2.
\]
Thus,
\[
\p_t \| \bu\|_{\dH^{m+1}}^2  \leq |\n \bu|_\infty Y_m +  \p^{m+1} \cC_{\a}(\bu,\rho)\cdot \p^{m+1}  \bu.
\]
Let us expand the commutator 
\[
\p^{m+1} \cC_{\a}(\bu,\rho) = \sum_{l=0}^{m+1} {m+1 \choose l} \cC_{\a} ( \p^l \bu, \p^{m+1-l} \rho).
\]
One end-point case, $l=m+1$, gives rise to a dissipative term:
\begin{multline*}
\int_{\T^n}  \cC_{\a}(  \p^{m+1} \bu,\rho)\cdot \p^{m+1}  \bu \dx = -  \frac12 \int_{\T^{2n}} \phi_{\a}(z) |\d_z \p^{m+1}  \bu(x)|^2 \rho(x+z) \dz \dx\\
- \frac12  \int_{\T^{2n}} \phi_{\a}(z) \d_z \p^{m+1} \bu(x) \p^{m+1} \bu(x)\d_z \rho(x) \dz \dx.
\end{multline*}
The first term is bounded by
\[ - \rmin \int_{\T^{2n}} \phi_{\a}(z)|\d_z \p^{m+1}  \bu(x)|^2   \dz \dx \sim - \rmin  \| \bu\|_{\dH^{m+1+ \frac{\a}{2}}}^2,
\]
which is the main dissipation term. The second is estimated as follows. Let us pick an $\e>0$ so small that $1+ \frac{\a}{2} > \a + \e$. Then
\[
\begin{split}
& \left| \int_{\T^{2n}} \phi_{\a}(z) \d_z \p^{m+1}  \bu(x)  \p^{m+1} \bu(x)\d_z \rho(x) \dz \dx \right| \leq  |\n \rho|_\infty \int_{\T^{2n}}  \frac{ |\p^{m+1} \d_z \bu(x)|}{|z|^{n/2 + \a - 1 + \e}} \frac{|\p^{m+1} \bu(x)|}{|z|^{n/2 - \e}} \dz \dx\\
& \leq |\n \rho|_\infty \|\bu\|_{H^{m+1}} \|\bu\|_{H^{m+\a +\e}} \leq  |\n \rho|_\infty \|\bu\|_{H^{m+1}} \|\bu\|_{H^{m+1+\a/2}} \leq \frac12 \rmin  \|\bu\|_{H^{m+1+\a/2}}^2 + \rmin^{-1} |\n \rho|_\infty^2 Y_m,
\end{split}
\]
where the first term is absorbed into dissipation. So,
\[
\int_{\T^n}  \cC_{\a}(  \p^{m+1} \bu,\rho)\cdot \p^{m+1}  \bu \dx \lesssim - \rmin  \| \bu\|_{\dH^{m+1+ \frac{\a}{2}}}^2 +  \rmin^{-1} |\n \rho|_\infty^2 Y_m.
\]
Let us consider first the other end-point case of $l=0$.  In this case the density suffers a derivative overload. We apply the following ``easing" technique:
\[
\int_{\T^n} \cC_{\a} ( \bu, \p^{m+1} \rho) \cdot \p^{m+1} \bu  \dx  = \int_{\T^{2n}} \phi_{\a}(z) \d_z \bu(x) \p^{m+1} \rho(x+z) \p^{m+1} \bu(x) \dz \dx.
\]
Observe that 
\[
\p^{m+1} \rho(x+z)  = \p_z \p_x^{m} \rho(x+z)  = \p_z ( \p_x^m \rho(x+z)  - \p_x^m \rho(x)) = \p_z \d_z \p^m \rho(x).
\]
Let us now integrate by parts in $z$:
\begin{multline*}
\int_{\T^n} \cC_{\a} ( \bu, \p^{m+1} \rho) \cdot \p^{m+1} \bu  \dx	=  \int_{\T^{2n}} \p_z \phi_{\a}(z) \d_z \bu(x)  \d_z \p^m \rho(x) \p^{m+1} \bu(x) \dz \dx \ + \\
+ \int_{\T^{2n}}  \phi_{\a}(z) \p \bu(x+z)   \d_z \p^m \rho(x) \p^{m+1} \bu(x) \dz \dx := J_1+ J_2.
\end{multline*}
Let us start with the  $J_2$ first. By symmetrization,
\[
\begin{split}
J_2 &=  \int_{\T^{2n}} \d_z \p \bu(x)   \d_z \p^m \rho(x) \p^{m+1} \bu(x) \phi_{\a}(z) \dz \dx -  \int_{\T^{2n}}  \p \bu(x)   \d_z \p^m \rho(x) \d_z \p^{m+1} \bu(x) \phi_{\a}(z) \dz \dx\\
&: = J_{2,1} + J_{2,2}.
\end{split}
\]
Term $J_{2,1}$ will appear in a series of similar terms that we will estimate systematically below. The bound for $J_{2,2}$ is rather elementary:
\[
J_{2,2} \leq  |\n \bu|_\infty \| \bu\|_{\dH^{m+1+ \a/2}} +  \| \rho\|_{H^{m+\a/2}} \leq \e \rmin  \| \bu\|_{\dH^{m+1+ \a/2}}^2 + \rmin^{-1} |\n \bu|_\infty^{2} Y_m.
\]
Similar computation can be made for $J_1$. Indeed, using that $\p_z \phi_{\a}(z)$ is odd, by symmetrization, we have
\[
J_1 = \frac12 \int_{\T^{2n}} \p_z \phi_{\a}(z) \d_z \bu(x)  \d_z \p^m \rho(x) \d_z \p^{m+1} \bu(x) \dz \dx.
\]
Replacing $|\d_z \bu(x) | \leq |z| | \bu |_\infty$, the rest of the term is estimated exactly as $J_{2,2}$. 

To summarize, we have obtained the bound
\[
\int_{\T^n} \cC_{\a} ( \bu, \p^{m+1} \rho) \cdot \p^{m+1} \bu  \dx \leq \e \rmin  \| \bu\|_{\dH^{m+1+ \a/2}}^2 + \rmin^{-1} |\n \bu|_\infty^{2} Y_m.
\]

Let us now examine the rest of the commutators $\cC_{\a} ( \p^l \bu, \p^{m+1-l} \rho)$ for $l=1,\ldots, m$. After symmetrization we obtain
\begin{multline*}
\int_{\T^n} \cC_{\a} ( \p^l \bu, \p^{m+1-l} \rho) \cdot \p^{m+1} \bu  \dx = \frac12 \int_{\T^{2n}} \d_z \p^l \bu(x)   \d_z \p^{m+1-l} \rho(x) \p^{m+1} \bu(x) \phi_{\a}(z) \dz \dx+\\
+ \frac12 \int_{\T^{2n}} \d_z \p^l \bu(x)   \p^{m+1-l} \rho(x) \d_z \p^{m+1} \bu(x) \phi_{\a}(z) \dz \dx: = J_1 + J_2.
\end{multline*}
Estimates on the new terms,  $J_1, J_2$ are a little more sophisticated as we seek to optimize distribution of $L^p$-norms inside their components. Notice that the case $l=1$ corresponds to the previously appeared term $J_{2,1}$. 

So, let us assume that $l=1,\ldots, m$.  We will distribute the parameters in $J_1$ as follows
\[
J_1 = \int_{\T^{2n}} \frac{\d_z \p^l \bu(x)}{|z|^{\frac{n}{p} + \frac{\a}{2} + 2\d }}   \frac{\d_z \p^{m+1-l} \rho(x)}{|z|^{\frac{n}{q} + \frac{\a}{2} }} \frac{ \p^{m+1} \bu(x)}{|z|^{\frac{n}{2} - \d}}  \frac{1}{|z|^{\frac{n}{r} - \d}}\dz \dx,
\]
where $\d >0$ is a small parameter to be determined later, and $(2,p,q,r)$ is a H\"older quadruple defined by
\[
p = 2 \frac{m + \frac{\a}{2}}{l-1 + \frac{\a}{2} }, \quad q = 2 \frac{m + \a-1}{m-l + \frac{\a}{2} }, \quad \frac{1}{r} = 1 - \frac{1}{2} - \frac{1}{p} - \frac{1}{q}.
\]
The existence of finite $r$ is warranted by the strict inequality which is verified directly:
\[
\frac{1}{2} + \frac{1}{p} + \frac{1}{q} < 1.
\] 
By the H\"older inequality,
\[
J_1 \leq  \|\bu \|_{\dW^{l+ \frac{\a}{2} + 2\d ,p }} \|\rho \|_{\dW^{ m+1 - l +\frac{\a}{2} ,q }} \|\bu \|_{\dH^{m+1}}.
\]
Let us apply the following Gagliardo-Nirenberg inequalities to all the terms
\[
\begin{split}
\|\bu \|_{\dH^{m+1}} & \leq \|\bu \|_{\dH^{m+1+ \frac{\a}{2}}}^{\frac{2m}{2m+\a}} | \n \bu |_2^{\frac{\a}{2m+\a}} \leq \|\bu \|_{\dH^{m+1+ \frac{\a}{2}}}^{\frac{2m}{2m+\a}} | \n \bu |_\infty^{\frac{\a}{2m+\a}} \\ 
\|\bu \|_{\dW^{l+ \frac{\a}{2} + 2\d ,p }} &\leq \|\bu \|_{H^{m+1 + \frac{\a}{2}}}^{\th_1} |\n \bu|_\infty^{1-\th_1}\\
\|\rho \|_{\dW^{ m+1 - l + \frac{\a}{2} ,q }} & \leq \|\rho \|_{H^{m+\a}}^{\th_2} |\n \rho|_\infty^{1-\th_2},
\end{split}
\]
where
\[
\th_1 = \frac{l-1 + \frac{\a}{2}  - \frac{n}{p} + 2\d}{m+\frac{\a}{2} - \frac{n}{2}}, \quad   \th_2 =  \frac{m-l+\frac{\a}{2}  - \frac{n}{q}}{m+\a-1 - \frac{n}{2}}.
\]
The exponents satisfy  the necessary requirements
\[
1\geq \th_1 \geq  \frac{l-1 + \frac{\a}{2}  + 2\d}{m+\frac{\a}{2} }, \quad 1 \geq  \th_2 = \frac{m-l+\frac{\a}{2}}{m+\a-1},
\]
and in fact,
\[
\th_1 = \frac{l-1 + \frac{\a}{2} }{m+\frac{\a}{2} } + O(\d).
\]
Now, we have
\[
J_1 \leq  \|\bu \|_{\dH^{m+1+ \frac{\a}{2}}}^{\frac{2m}{2m+\a} + \th_1}    \|\rho \|_{H^{m+\a}}^{\th_2}     | \n \bu |_\infty^{\frac{\a}{2m+\a} + 1 - \th_1}   |\n \rho|_\infty^{1-\th_2}.
\]
By generalized Young,
\[
J_1 \leq  \e \rmin  \| \bu\|_{\dH^{m+1+ \a/2}}^2  + \rmin^{-1}   \|\rho \|_{H^{m+\a}}^{\th_2 Q}  (| \n \bu |_\infty^{\frac{\a}{2m+\a} + 1 - \th_1}   |\n \rho|_\infty^{1-\th_2})^Q,
\]
where  $Q$ is the conjugate to  $\frac{2m}{2m+\a} + \th_1$. We have $\th_2 Q< 2$ as long as 
\begin{equation*}
\th_1 + \th_2 < 2 - \frac{2m}{2m+\a}.
\end{equation*}
We in fact have even stronger inequality, $\th_1 + \th_2 <1$ provided $\d$ is small enough. So, we arrived at
\[
J_1 \leq \e \rmin  \| \bu\|_{\dH^{m+1+ \a/2}}^2  + \rmin^{-1} p_N(|\n \rho|_\infty,| \n \bu |_\infty) Y_m,
\]
for some polynomial $p_N$. 

Finally, moving on to $J_2$, we distribute the exponents as follows 
\[
J_2 \leq \int_{\T^{2n}} \frac{|\d_z \p^l \bu(x)|}{|z|^{\frac{n}{p} + 2\d + \frac{\a}{2}}}   \frac{| \p^{m+1-l} \rho(x)| }{|z|^{\frac{n}{q} - \d}} \frac{ |\d_z \p^{m+1} \bu(x)|}{ |z|^{\frac{n}{2} + \frac{\a}{2} }} \frac{1}{|z|^{\frac{n}{r} -\d}}\dz \dx \leq  \|\bu \|_{\dW^{l+\d+ \frac{\a}{2},p }} \|\rho \|_{\dW^{ m+1 - l,q }} \|\bu\|_{\dH^{m+1+\frac{\a}{2}}} .
\]
Here we choose $(r,p,q,\d)$ as follows
\[
q = 2 \frac{m+\a-1}{m-l}, \quad p = 2 \frac{m + \frac{\a}{2}}{l-1 + \frac{\a}{2}},\quad \frac{1}{r} = 1 - \frac{1}{2} - \frac{1}{p} - \frac{1}{q},
\]
and  $\d$ is small.  With these choices we proceed with the  Gagliardo-Nirenberg inequalities
\[
\begin{split}
\|\bu \|_{\dW^{l+2\d+ \frac{\a}{2},p }}  & \leq \|\bu\|_{\dH^{m+1+\frac{\a}{2}}}^{\th_1} |\n \bu|_\infty^{1-\th_1} \\
\|\rho \|_{\dW^{ m+1 - l,q }}  & \leq \|\rho \|_{\dH^{ m+\a }}^{\th_2} |\n \rho|_\infty^{1-\th_2},
\end{split}
\]
where
\[
\th_1 = \frac{l-1 + \frac{\a}{2} + 2\d}{m + \frac{\a}{2} - \frac{n}{2} }=  \frac{l-1 + \frac{\a}{2} }{m+\frac{\a}{2} } + O(\d), \quad \th_2 = \frac{m-l}{m+\a-1}.
\]
Now to achieve the bound
\[
J_2 \leq \e \rmin  \| \bu\|_{\dH^{m+1+ \a/2}}^2  + \rmin^{-1} p_N(|\n \rho|_\infty,| \n \bu |_\infty) Y_m,
\]
we have to make sure that  $\th_1 + \th_2 \leq 1$, which is true for small $\d$.

We have proved the following a priori bound on $\bu$:
\[
\p_t \| \bu\|_{\dH^{m+1}}^2  \leq - \frac12 \rmin  \| \bu\|_{\dH^{m+1+ \a/2}}^2 + \rmin^{-1} p_N(|\n \rho|_\infty,| \n \bu |_\infty) Y_m.
\]
Together with the previously established bounds we obtain
\[
\ddt Y_m \leq - \frac12 \rmin  \| \bu\|_{\dH^{m+1+ \a/2}}^2 + \rmin^{-1} p_N(|\n \rho|_\infty,| \n \bu |_\infty, |e|_\infty) Y_m.
\]
This of course implies a Riccati inequality, provided $m > \frac{n}{2}+1$:
\[
\ddt Y_m \leq C Y_m^N,
\]
and provides a priori bound independent of the viscosity coefficient. Thus, we can extend it to an interval  independent of $\e$ as well. By the compactness argument similar to the smooth kernel case, we obtain a local solution in the same class as initial data and $\bu \in L^2([0,T_0); \dH^{m+1+ \a/2})$. In addition, we obtain a continuation criterion  -- as long as  $|\n \rho|_\infty,| \n \bu |_\infty$, $|e|_\infty$ remain bounded on $[0,T_0)$ the solutions can be extended beyond $T_0$.  However  everything is reduced to a control over the first two quantities, because $|e|_\infty$ remains bounded as long as  $|\n \bu|_\infty$ is in view of \eqref{e:euu}. 

It is clear from the proof that \eqref{e:BKMsing} can be replaced with an integrability condition with some high power depending on $m,n,\a$.

\section{Global well-posedness and strong flocking}\label{s:GWP}
According to our local well-posedness Theorem \ref{t:lwpsing} we already have a local solution $(u,\rho)$ on time
interval $[0, T_0)$. We proceed in several steps. First, we establish uniform bounds \eqref{e:rhopositive} on the density which
depend only on the initial conditions. So, such bounds hold uniformly on the available time interval $[0, T_0)$.
Next, we invoke results from the theory fractional parabolic equations to conclude that our solution gains H\"{o}lder regularity after a short period of time, and the H\"{o}lder exponent as well as the bound on the H\"{o}lder
norm depend on the $L^{\infty}$ bound of the solution. Finally, we establish a continuation criterion much weaker
than that of Theorem \ref{t:lwpsing} -- claiming that any H\"{o}lder regularity of the density propels higher order norms beyond $T_0$. \\

Paired with the density equation we find that the ratio $q:=e/\rho$ satisfies the transport equation
\begin{equation}\label{e:Dtq}
\DDt q = q_t + u q_1=0.
\end{equation}
Starting from sufficiently smooth initial condition with $\rho_0$ away from vacuum we can assume that
\begin{equation}\label{e:q(t)=q(0)}
|q(t)|_{\infty}=|q_0|_{\infty}<\infty.
\end{equation}
\noindent {\sc Step 1: bounds on the density.}  
We start by establishing uniform bounds \eqref{e:rhopositive} on the density which depend only on the initial conditions.
First, recall that $q$ is transported \eqref{e:Dtq}, and hence is bounded for all time with its initial
value $|q_0|_{\infty}$. So, we can write the continuity equation  as
\[
\rho_t + u \rho_1 =-q \rho^2-\rho \L_{\a}(\rho).
\]
Let us evaluate at a point $x^{+}(t)$ where the maximum of $\rho(\cdot,t)$, denoted by $\rho^{+}(t):=\rho(x^{+}(t),t)$, is reached. We obtain
\begin{align*}
\ddt \rho^{+}(t)&=-q(x^{+}(t),t)\left(\rho^{+}(t)\right)^2-\rho^{+}(t)\int_{\R^n} \left(\rho^{+}(t)- \rho(x^{+}(t)+z,t)\right)\frac{\dz}{|z|^{n+\a}}\\
&\leq |q_0|_{\infty}\left(\rho^{+}(t)\right)^2-\frac{\rho^{+}(t)}{r^{n+\a}}\int_{|z|<r} \left(\rho^{+}(t)- \rho(x^{+}(t)+z,t)\right)\dz\\
&\leq |q_0|_{\infty}\left(\rho^{+}(t)\right)^2-\frac{\rho^{+}(t)}{r^{n+\a}}\left(V_n(r)\rho^{+}(t)-\cM \right),
\end{align*}
where $V_n(r)$ denotes the $n$-dimensional volume of a ball of radius $r$. As $V_n(r)=C(n)r^n$, we get
\[
\ddt \rho^{+}(t)\leq \left[|q_0|_{\infty}-\frac{C(n)}{r^{\a}}\right]\left(\rho^{+}(t)\right)^2+\frac{\cM}{r^{n+\a}}\rho^{+}(t).
\]
Let us pick $r$ small enough so that $\frac{C(n)}{r^{\a}}\geq |q_0|_{\infty}+1$. Then 
\[
\ddt \rho^{+}(t)\leq - \left(\rho^{+}(t)\right)^2 +c_0\rho^{+}(t)
\] 
which establishes the upper bound by integration.

As to the lower bound we argue similarly. Let $\rho^{-}(t)$ be the minimum value of $\rho(\cdot,t)$ and $x^{-}(t)$ a point
where such value is achieved. We have 
\begin{align*}
\ddt \rho^{-}(t)&=-q(x^{-}(t),t)\left(\rho^{-}(t)\right)^2-\rho^{-}(t)\int_{\T^n} \left(\rho^{-}(t)- \rho(x^{-}(t)+z,t)\right)\phi_{\a}(z)\dz\\
&\geq -|q_0|_{\infty}\left(\rho^{-}(t)\right)^2-\phi_{\a}^{-}\,\rho^{-}(t)\int_{\T^n} \left(\rho^{-}(t)- \rho(x^{-}(t)+z,t)\right)\dz\\
&\geq -|q_0|_{\infty}\left(\rho^{-}(t)\right)^2-\phi_{\a}^{-}\,\rho^{-}(t)\left((2\pi)^n\rho^{-}(t)-\cM\right).
\end{align*}
Note that at this point the global communication of the model is crucial: $\phi_{\a}^{-}:=\inf_{z\in\T^n}\phi_{\a}(z)>0.$
Then
\[
\ddt \rho^{-}(t)\geq -c_1 \left(\rho^{-}(t)\right)^2+c_2 \rho^{-}(t)
\]
which establishes the lower bound by integration.\\

\noindent {\sc Step 2: bounds on the entropy.} 
As an immediate consequence of the uniform bound on the density and \eqref{e:q(t)=q(0)} we have a uniform global bound on the entropy $|e(t)|_{\infty}<\infty$. This argument can be iterated to higher derivatives as follows.
Let us start with one observation in the $x_1$-direction. Note that if a quantity $q$ is transported \eqref{e:Dtq}, then the same transport equation governs $q_1/\rho$
\begin{equation*}
\DDt \left(\frac{q_1}{\rho} \right) = \left(\frac{q_1}{\rho} \right)_t + u \left(\frac{q_1}{\rho} \right)_1=0.
\end{equation*}
Consequently, if $|q_1|/\rho$ is bounded at initial time $t=0$ it will remain bounded at later time $t>0$. Unraveling the formulas, we obtain the bound
\begin{equation*}
|e_1(t)|_{\infty}\lesssim |\rho_1(t)|_{\infty}.
\end{equation*}
For the rest of derivatives of the entropy $e_i$ with $i=2,\ldots,n$, we have that
\begin{equation}\label{e:Dtqi/rho}
\DDt \left(\frac{q_i}{\rho} \right) = \frac{\rho \mbox{D}_t q_i -q_i \mbox{D}_t\rho }{\rho^2}=\left(\frac{q_i}{\rho} \right)u_1-\left(\frac{q_1}{\rho} \right)u_i.
\end{equation}
The proof is just a combination of the following two facts.
On one hand, as $q$ satisfies the transport equation \eqref{e:Dtq}, the first material derivative   reduces to $\mbox{D}_t q_i \equiv (\mbox{D}_t q)_i-u_i q_1=-u_i q_1$. 
On the other hand, the second material derivative  is just $\mbox{D}_t \rho \equiv\rho_t +(u\rho)_1 -u_1 \rho = -u_1 \rho$ thanks to the continuity equation of the density. Consequently, applying Gr\"{o}nwall's inequality in \eqref{e:Dtqi/rho} give us the bound
\begin{equation}\label{e:grade}
|\nabla e (t)|_{\infty} \lesssim |\nabla \rho (t)|_{\infty}+ |\nabla u|_{L^1_{t} L^{\infty}_{x}},
\end{equation}
using that $q_1/\rho$ is preseved in time. Let us note that in order to make pointwise evaluation possible in \eqref{e:grade} one has to assume regularity $\nabla e\in H^{m-1}(\T^n)\subseteq C(\T^n)$ which guaranteed provided \nolinebreak $m>\tfrac{n}{2}+1.$\\

\noindent {\sc Step 3: H\"{o}lder regularization.} 
The parabolic nature of the density equation is an essential structural feature of the system that has been used in all of the preceding works in 1D. Using the $e$-quantity we can write
\begin{equation}\label{e:density}
\rho_t + u \rho_1 +e \rho=-\rho \L_{\a}(\rho).
\end{equation}
Similarly, one can write the equation for the momentum $m=\rho u$:
\begin{equation}\label{e:momemtum}
m_t+u m_1 +e m=-\rho\Lambda_{\a}(m).
\end{equation}
Note that in both cases the drift $u$ and the forcing $e\rho$  or $em$ are bounded a priori due to the maximum principle stated above.  Hence, the density and momemtum equations \eqref{e:density}, \eqref{e:momemtum} falls under the general class of  forced fractional parabolic equations with bounded drift and force:
\[
w_t+u\cdot\nabla w=\mathcal{L}_{\a}(w)+f, \qquad \mathcal{L}_{\a}(w)(x,t):=\int_{\R^n}\mbox{K}(x,z,t)\left(w(x+z,t)-w(x,t)\right) \dz
\]
with a diffusion operator associated with the singular kernel $\mbox{K}(x,z,t)=\rho(x,t)|z|^{-(n+\a)}$ which is even with respect to $z$. The bounds on the density provide uniform ellipticity bounds on
the kernel:
\[
(2-\a)\frac{\lambda}{|z|^{n+\a}}\leq \mbox{K}(x,z,t)\leq (2-\a)\frac{\Lambda}{|z|^{n+\a}}
\]
The most common assumption in the literature is that for all $x$ and $t$, the kernel $\mbox{K}$ is comparable pointwise in terms of $z$ to the kernel for the fractional Laplacian.

Regularity of these equations has been the subject of active research in recent years. In particular, the result of Silvestre \cite{S2012}, see also Schwab and Silvestre \cite{SS2016},  which provides H\"{o}lder regularization bound for some $\gamma>0$ given by
\begin{align}\label{e:gamma-rho}
\|\rho\|_{C^{\gamma}(\T^n \times [T/2,T))} &\lesssim |\rho|_{L^{\infty}(\T^n \times [0,T))}+|\rho e|_{L^{\infty}(\T^n \times [0,T))},\\
\|m\|_{C^{\gamma}(\T^n \times [T/2,T))} &\lesssim |m|_{L^{\infty}(\T^n \times [0,T))}+|m e|_{L^{\infty}(\T^n \times [0,T))}, \nonumber
\end{align}
and
\begin{equation}\label{e:gamma-u}
\|u\|_{C^{\gamma}(\T^n \times [T/2,T))} \leq C\left(|\rho|_{L^{\infty}(\T^n \times [0,T))},|u|_{L^{\infty}(\T^n \times [0,T))}\right),
\end{equation}
where the latter inequality follows from \eqref{e:gamma-rho} since $\rho$ is bounded below.
Since the right hand side  of \eqref{e:gamma-rho} and \eqref{e:gamma-u} is uniformly bounded on time we have obtained uniform bound on $C^{\gamma}$-norm starting, by rescaling, from any positive time. These results apply in our case when $1\leq \a<2$ since in this case we only need bounded drift and force. Due to futher limitations that will come later we will only proceed with $1<\a<2$, however, this initial regularization technically holds even for $\a=1$.

We now proceed to establishing that the solution fulfills the above continuation criterion \eqref{e:BKMsing}:  $|\nabla \rho|_{\infty}+|\nabla u|_{\infty} \in L^\infty([0,T_0))$, where $[0,T_0)$ is a given local interval of existence.

\noindent {\sc Step 4.1: Control over $|\nabla \rho|_{\infty}$.} So, let us start with $\p \rho \equiv \rho_{i}$ for $i\in \lbrace 1,\ldots,n\rbrace.$ Then, we have
\[
\p_t \p \rho + \p \p_1 ( u\rho )=0,
\]
or expanding the non-linear part in higher and lower order terms, we arrive at
\[
\p_t  \rho_i + (\rho  u_{1i} + u \rho_{1i}) + (u_i \rho_1 + u_1 \rho_i)=0.
\]
Evaluating at the maximum of $|\rho_{i}|$ and multiplying by $\rho_{i}$ again (we use the classical Rademacher theorem here to justify the time derivative) we obtain
\[
\p_t |\rho_i|^2 + \rho \rho_i u_{1i}+(u_i \rho_1 + u_1 \rho_i)\rho_i =0.
\]
Using the entropy $e =  u_1 -\Lambda_{\a}\rho$, we write the remaining higher order term in a more convenient \nolinebreak way.
\[
\p_t |\rho_i|^2 + \rho \rho_i e_{i}+(u_i \rho_1 + u_1 \rho_i)\rho_i =-\rho \rho_i \Lambda_{\a}(\rho_i).
\]
In consequence, summing over indexes $i\in\lbrace 1,\ldots,n\rbrace$ and using the upper and lower uniform bounds previously proved  for the density, we obtain that the following estimate holds:
\[
\p_t |\nabla \rho |_{\infty}^2 \leq \rho^{+}(t)|\nabla \rho|_{\infty} |\nabla e|_{\infty} +2 |\nabla u|_{\infty} |\nabla \rho|_{\infty}^2-\rho^{-}(t)\nabla \rho(x^{\star}) \cdot \Lambda_{\a}(\nabla \rho)(x^{\star}).
\]
Next, in view of the pointwise identity \cite{CC2004}, we have
\[
\nabla f(x)\cdot \Lambda_{\a}(\nabla f)(x)= \tfrac{1}{2}\Lambda_{\a}(|\nabla f|^2)(x)+\tfrac{1}{2} D_{\a}(\nabla f)(x)
\]
where
\[
D_{\a}(\nabla f)(x):=\int_{\R^n}\frac{|\nabla f(x+z)-\nabla f(x)|^2}{|z|^{n+\a}} \mbox{d}z.
\]
In additon, using the non-linear bounds from \cite{CV2012} the following pointwise bound holds
\begin{equation}\label{e:NLmax}
D_{\a}(\nabla f)(x) \gtrsim \frac{|\nabla f(x)|^{2+\a}}{ | f|_{\infty}^{\a}}.
\end{equation}
The above non-local maximum principle yields the following bound at the maximal point $x^{\star}\equiv x^{\star}(t)$:
\[
\nabla \rho(x^{\star}) \cdot \Lambda_{\a}(\nabla \rho)(x^{\star})\geq \frac{1}{4} D_{\a}(\nabla \rho)(x^{\star})+\frac{c}{|\rho|_{\infty}^{\a}} |\nabla \rho|_{\infty}^{2+\a}.
\]
Due to the uniform bound from below on $\rho$, we arrive at
\begin{equation}\label{e:finalgradrho}
\p_t |\nabla \rho |_{\infty}^2 \lesssim C \left(|\nabla \rho|_{\infty} |\nabla e|_{\infty} + |\nabla u|_{\infty} |\nabla \rho|_{\infty}^2\right)- D_{\a}(\nabla \rho)(x^{\star})- |\nabla \rho|_{\infty}^{2+\a}.
\end{equation}

\noindent {\sc Step 4.2: Control over $|\nabla u |_{\infty}$.} We continue with $\p u \equiv  u_{i}$ for $i \in \lbrace 1,\ldots,n\rbrace.$
\[
\p_t\p u +\p (u \p_1 u)=\p \cC_{\a}(u,\rho)
\]
where
\[
\cC_{\a}(f,g)=f\Lambda_{\a}(g)-\Lambda_{\a}(fg)=\text{p.v.}\int_{\R^n}\frac{\delta_{z}f(x)}{|z|^{n+\a}} g(x+z) \dz.
\]
We can rewrite the above expression as
\[
\p_t u_i+u_{i} u_1+u u_{1i}=\cC_{\a}(u_{i},\rho)+\cC_{\a}(u,\rho_{i}).
\]
Evaluating at the maximum of $|u_{i}|$ and multiplying by $u_{i}$ again (we use the classical Rademacher theorem here to justify the time derivative) we obtain
\[
\p_t| u_{i}|^2 +  |u_i|^2 u_1 = u_{i}\cC_{\a}(u_{i},\rho) + u_{i}\cC_{\a}(u,\rho_{i}).
\]
In consequence, summing over indexes $i\in\lbrace 1,\ldots,n\rbrace ,$ we obtain that the following estimate holds:
\[
\p_t |\nabla u|_{\infty}^2 \leq |\nabla u|_{\infty}^3 +\nabla u(x_{\star}) \cdot \cC_{\a}(\nabla u,\rho)(x_{\star})+ \nabla u (x_{\star}) \cdot \cC_{\a}( u,\nabla \rho)(x_{\star}).
\]
The dissipation term is bounded, as before by
\[
\nabla u (x_{\star}) \cdot \cC_{\a}(\nabla u,\rho)(x_{\star})=\nabla u (x_{\star}) \cdot \int_{\R^n}\frac{\delta_{z}\nabla u(x_{\star})}{|z|^{n+\a}} \rho(x_{\star}+z)\dz \leq -\rho^{-}(t) D_{\a}(\nabla u)(x_{\star}).
\]
In addition, it obeys another non-local maximum principle similar to \eqref{e:NLmax} where instead of $\nabla u$ we replace it with $\nabla (u-\bar{u})$, thus the denominator contain the amplitude $\cA(t)$ rather than $|u|$:
\[
D_{\a}(\nabla u)(x_{\star}) \gtrsim  \frac{|\nabla u|_{\infty}^{2+\a}}{ \cA^{\a}(t)}.
\]
As before, due to the uniform bound from below on $\rho$, we arrive at
\begin{equation}\label{e:graduaux}
\p_t |\nabla u|_{\infty}^2 \lesssim C\left(|\nabla u|_{\infty}^3+  \nabla u (x_{\star}) \cdot \cC_{\a}( u,\nabla \rho)(x_{\star})\right) - D_{\a}(\nabla u)(x_{\star})-\frac{|\nabla u|_{\infty}^{2+\a}}{\cA^{\a}(t)}.
\end{equation}
Note that the remaining higher order term can be write as
\[
\nabla u\cdot \cC_{\a}(u,\nabla\rho)=\nabla u \cdot \int_{\R^n} \frac{\delta_{z}u(x)}{|z|^{n+\a}} \nabla \rho(x+z)\dz.
\]
To handle it, we split the integral into two parts:
\begin{equation}\label{e:split}
\int_{\R^n} \frac{\delta_{z}u(x)}{|z|^{n+\a}} \nabla \rho(x+z) \dz= \int_{|z|<1} \frac{\delta_{z}u(x)}{|z|^{n+\a}} \nabla \rho(x+z) \dz+ \int_{|z|>1} \frac{\delta_{z}u(x)}{|z|^{n+\a}} \nabla \rho(x+z) \dz,
\end{equation}
where we use the alignment in the large scale part
\begin{equation}\label{e:large}
\int_{|z|>1} \frac{\delta_{z}u(x)}{|z|^{n+\a}} \nabla \rho(x+z) \dz \lesssim  \cA(t) |\nabla \rho|_{\infty}.
\end{equation}
To handle the small scale part, we fix a scale parameter $1 > r > 0$ to be determined later, and split the integral  into two different parts: small-small scale and small-middle scale. 
\[
\int_{|z|<1} \frac{\delta_{z}u(x)}{|z|^{n+\a}} \nabla \rho(x+z) \dz=\int_{|z|<r} \frac{\delta_{z}u(x)}{|z|^{n+\a}} \nabla \rho(x+z) \dz+\int_{r<|z|<1} \frac{\delta_{z}u(x)}{|z|^{n+\a}} \nabla \rho(x+z) \dz.
\]
If we add and subtract $\nabla u(x)\cdot z$ in the small-small scale part, we obtain
\begin{align*}
\int_{|z|<1} \frac{\delta_{z}u(x)}{|z|^{n+\a}} \nabla \rho(x+z) \dz&=\int_{|z|<r} \frac{[\delta_{z}u(x)-\nabla u(x)\cdot z]}{|z|^{n+\a}} \nabla \rho(x+z) \dz+\int_{r<|z|<1} \frac{\delta_{z}u(x)}{|z|^{n+\a}} \nabla \rho(x+z) \dz\\
&\quad + \nabla u(x)\cdot \int_{|z|<r} \frac{z}{|z|^{n+\a}} \nabla \delta_{z}\rho(x) \dz \equiv \cI_1 + \cI_2 +\cI_3.
\end{align*}
Then, we will study each term separately:\\
\underline{\sc $\cI_1$-Term:}  The most singular part of \eqref{e:split} has to be handled in a way that utilizes dissipation. First, let us use the analogue of spherical coordinates in $n$-dimension to write the dissipation term in the following form:
\begin{equation}\label{e:radialD}
\quad D_{\a}(f)(x)=\int_{\R^n}|\delta_{z}f(x)|^2\frac{\dz}{|z|^{n+\a}}=\int_{S^{n-1}}\int_{0}^{\infty} |\delta_{r \theta}f(x)|^2\frac{r^{n-1}\dr}{|r|^{n+\a}}\mbox{d}_{S^{n-1}} V
\end{equation}
where $\mbox{d}_{S^{n-1}} V = \sin^{n-2}(\theta_1)\sin^{n-3}(\theta_2)\ldots \sin(\theta_{n-2}) \mbox{d} \theta_1 \mbox{d} \theta_2 \ldots \mbox{d} \theta_{n-1}$ is the volume element of  $S^{n-1}$, which is the generalization of the ordinary sphere to spaces of arbitrary dimension. In the rest, we will use the notation
\[
D_{\a}(f)(x,\theta):=\int_{0}^{\infty} |\delta_{r \theta}f(x)|^2\frac{\dr}{|r|^{1+\a}}
\]
to denote the inner radial integral of \eqref{e:radialD}. Then,  the dissipation can be write as
\[
D_{\a}(f)(x)=\int_{S^{n-1}}D_{\a}(f)(x,\theta) \mbox{d}_{S^{n-1}} V.
\]
Let us estimate a Taylor expansion in terms of $D_{\a}(\nabla f)(x)$. First, using the identity
\[
\delta_{z}f(x)=\int_{0}^1 \nabla f(x+sz)\cdot z \ds,
\]
we have that
\[
\delta_{z}f(x)-\nabla f(x)\cdot z=\int_0^1 \left[\nabla f (x+s z)-\nabla f(x)\right]\cdot z\, \ds= \int_{0}^{z} \left[\nabla f (x+ w)-\nabla f(x)\right]\cdot  \dw
\]
where the last integral is over the straight radial segment $[0,z]\subset \R^n$. Then, after make the change of variables, we obtain for $\bar{\theta}\equiv\bar{\theta}(z)\in S^{n-1}$ that
\[
\delta_{z}f(x)-\nabla f(x)\cdot z = \int_{0}^{z} \delta_{w}\nabla f (x) \cdot \dw=\int_{0}^{|z|}\delta_{r\bar{\theta}}\nabla f (x) \dr
\]
and H\"{o}lder inequality give us that
\begin{equation*}
|\delta_{z}f(x)-\nabla f(x)\cdot z|\leq \sqrt{D_{\a}(\nabla f)(x,\bar{\theta})} |z|^{(2+\a)/2}.
\end{equation*}
After this, the small-small scale part can be handled using dissipation as follows
\begin{equation}\label{e:I1term}
\cI_1\equiv \int_{|z|<r}\frac{[\delta_{z}u(x)-\nabla u(x)\cdot z]}{|z|^{n+\a}}\nabla \rho(x+z) \dz \lesssim |\nabla \rho|_{\infty} \sqrt{D_{\a}(\nabla u)(x)} r^{1-\tfrac{\a}{2}}.
\end{equation}

\noindent
\underline{\sc $\cI_2$-term:}  In the small-middle scale part we use the available H\"{o}lder regularity for $\a\in (1,2).$
\begin{equation}\label{e:I2term}
\cI_2\equiv  \int_{r<|z|<1}\frac{\delta_{z}u(x)}{|z|^{n+\a}} \nabla \rho(x+z)\dz  \lesssim |\nabla \rho|_{\infty} \|u\|_{C^{\gamma}} r^{\gamma-\a}.
\end{equation}

\noindent
\underline{\sc $\cI_3$-term:} Finally, the remaining term can be rewrite as
\begin{align*}
\nabla u(x)\cdot \int_{|z|<r} \frac{z}{|z|^{n+\a}} \nabla \delta_{z} \rho(x) \dz&= \nabla u(x)\cdot \left[\int_{\R^n} \frac{z}{|z|^{n+\a}} \nabla \delta_{z}\rho(x) \dz- \int_{|z|>r} \frac{z}{|z|^{n+\a}} \nabla \delta_{z} \rho(x) \dz\right]\\
&= \nabla u(x) \Lambda_{\a}(\rho)(x)-\nabla u(x)\cdot \int_{|z|>r} \frac{z}{|z|^{n+\a}} \nabla \delta_{z} \rho(x) \dz
\end{align*}
where the full integral over $\R^n$ is nothing other (by  integration by parts and the fact that $|z|^{-(n+\a)}\approx\nabla\cdot(z |z|^{-(n+\a)})$) than the integral representation of the fractional Laplacian.

In order to handle the remaining residual term, we introduce the usual even cut-off function $\psi\in C^{\infty}$ with $\psi(z)=1$ for $|z|<1$ and $\psi(z)=0$ for $|z|>2$. Denote $\psi_r(z)=\psi(z/r)$, and decompose
\[
\int_{|z|>r} \frac{z}{|z|^{n+\a}} \nabla \delta_{z} \rho(x) \dz=\int_{\R^n} (1-\psi_r(z))\frac{z}{|z|^{n+\a}} \nabla \delta_{z}\rho(x)\dz+\int_{r<|z|<2r}\psi_r(z) \frac{z}{|z|^{n+\a}} \nabla \delta_{z} \rho(x) \dz. 
\]
After apply integration by parts and use the  available H\"{o}lder regularity, we obtain
\[
\int_{|z|>r} \frac{z}{|z|^{n+\a}} \nabla \delta_{z} \rho(x) \dz\lesssim |\nabla \rho|_{\infty}+\|\rho\|_{C^{\gamma}} r^{\gamma-\a}. 
\]
Finally, by the uniform $L^{\infty}$  bound previously proved for the entropy and the expression $e=u_1-\Lambda_{\a}\rho$ we arrive at
\begin{equation}\label{e:I3term}
\cI_3\equiv \nabla u(x)\cdot \int_{|z|<r} \frac{z}{|z|^{n+\a}} \nabla \delta_{z} \rho(x) \dz \lesssim |\nabla u|_{\infty}^2+|\nabla u|_{\infty} |\nabla \rho|_{\infty} +|\nabla u|_{\infty}\|\rho\|_{C^{\gamma}} r^{\gamma-\a}.
\end{equation}

Combining \eqref{e:I1term}, \eqref{e:I2term} and \eqref{e:I3term}  we have proved that the following estimate holds:
\begin{align}\label{e:small}
\int_{|z|<1} \frac{\delta_{z}u(x)}{|z|^{n+\a}} \nabla \rho(x+z) \dz &\lesssim |\nabla \rho|_{\infty} \sqrt{D_{\a}(\nabla u)(x)} r^{1-\tfrac{\a}{2}}+\left(|\nabla \rho|_{\infty} +|\nabla u|_{\infty}\right) r^{\gamma-\a}\\
&\quad +|\nabla u|_{\infty}\left(|\nabla u|_{\infty} +|\nabla \rho|_{\infty}\right). \nonumber
\end{align}
The competition occurs only between the small and middle range terms. Optimizing over r we set
\[ 
r\approx \left( D_{\a}(\nabla u)(x)\right)^{-\frac{1}{2+\a-2\gamma}}
\]
unless such expression is $>1$, in which case we have an absolute bound on the dissipation and the proof proceeds trivially. With the established bounds we obtain the following pointwise estimate
\begin{align*}
\int_{|z|<1} \frac{\delta_{z}u(x)}{|z|^{n+\a}} \nabla \rho(x+z) \dz &\lesssim \left(|\nabla \rho|_{\infty} +|\nabla u|_{\infty}\right) \left[\left(  D_{\a}(\nabla u)(x)\right)^{\frac{\a-\gamma}{2+\a-2\gamma}}+|\nabla u|_{\infty}\right].
\end{align*}
Consequently, combining \eqref{e:small} and \eqref{e:large} into \eqref{e:split} we have proved the following estimate:
\begin{align}\label{e:eatD+cubic}
\nabla u (x) \cdot \cC_{\a}( u,\nabla \rho)(x)&\lesssim \left(|\nabla u|_{\infty}^{2}+ |\nabla \rho|_{\infty}^{2}\right)\left(  D_{\a}(\nabla u)(x)\right)^{\frac{\a-\gamma}{2+\a-2\gamma}}\\
&\quad +|\nabla u|_{\infty}^2\left(|\nabla u|_{\infty} +|\nabla \rho|_{\infty}\right). \nonumber
\end{align}
Note that for $\a\in (0,2)$ and $\gamma>0$  we have that $\frac{\a-\gamma}{2+\a-2\gamma} < \frac{\a}{2+\a}.$ So, we can use the generalized Young inequality to obtain
\begin{equation}\label{e:eatD}
\left(|\nabla u|_{\infty}^{2}+ |\nabla \rho|_{\infty}^{2}\right)\left(  D_{\a}(\nabla u)(x)\right)^{\frac{\a-\gamma}{2+\a-2\gamma}}\lesssim c_{\e}+ \e\left(|\nabla \rho|_{\infty}^{2+\a}+\frac{|\nabla u|_{\infty}^{2+\a}}{\cA^{\a}(t)}+ D_{\a}(\nabla u)(x)\right).
\end{equation}
At this point, adding \eqref{e:eatD} to \eqref{e:eatD+cubic}, we obtain 
\begin{align*}
\nabla u (x_{\star}) \cdot \cC_{\a}( u,\nabla \rho)(x_{\star})&\leq  C\, |\nabla u|_{\infty}^2\left( |\nabla \rho|_{\infty}+|\nabla u|_{\infty}\right) +c_{\e}\\
&\quad + \e\left(|\nabla \rho|_{\infty}^{2+\a}+\frac{|\nabla u|_{\infty}^{2+\a}}{\cA^{\a}(t)}+ D_{\a}(\nabla u)(x_{\star})\right)
\end{align*}
and plugging this into \eqref{e:graduaux} we arrive at
\begin{align}\label{e:finalgradu}
\p_t |\nabla u|_{\infty}^2 &\lesssim C\, |\nabla u|_{\infty}^2\left( |\nabla \rho|_{\infty}+|\nabla u|_{\infty}\right) +\e |\nabla \rho|_{\infty}^{2+\a} + c_{\e}\\
&\quad - (1-\e)D_{\a}(\nabla u)(x_{\star})-(1-\e)\frac{|\nabla u|_{\infty}^{2+\a}}{\cA^{\a}(t)}. \nonumber
\end{align}

\noindent {\sc Step 4.3: Control over $|\nabla \rho|_{\infty}+|\nabla u|_{\infty}$.} In order to complete the proof of Theorem \ref{t:flock-singular}, we combine \eqref{e:finalgradrho} and \eqref{e:finalgradu} to obtain  that
\begin{align}\label{e:grad(rho+u)}
\p_t \left[|\nabla \rho |_{\infty}^2+|\nabla u |_{\infty}^2\right] &\lesssim C \left(|\nabla \rho|_{\infty} |\nabla e|_{\infty} +|\nabla u|_{\infty}^{2} |\nabla \rho|_{\infty}+ |\nabla u|_{\infty} |\nabla \rho|_{\infty}^2 +|\nabla u|_{\infty}^3\right)+c_{\e}\\
& \quad  -(1-\e)\frac{|\nabla u|_{\infty}^{2+\a}}{\cA^{\a}(t)}-(1-\e)|\nabla \rho|_{\infty}^{2+\a} \nonumber \\
& \quad - (1-\e)D_{\a}(\nabla u)(x_{\star}) -D_{\a}(\nabla \rho)(x^{\star}). \nonumber
\end{align}
The above expression emphasizes the fact that the point at which each $D$-term is evaluated is different. Notice that for $\alpha>1$, we can absorb the cubic terms in \eqref{e:grad(rho+u)} simply by interpolation:
\begin{align}\label{e:cubic}
|\nabla u|_{\infty}^2 |\nabla \rho|_{\infty}&\leq \e \left( \frac{|\nabla u|_{\infty}^{2+\a}}{\cA^{\a}(t)}+|\nabla \rho|_{\infty}^{2+\a}\right)+c_{\e}\cA^{\frac{2\a}{\a-1}}(t),\nonumber\\
|\nabla u|_{\infty} |\nabla \rho|_{\infty}^2 &\leq \e \left( \frac{|\nabla u|_{\infty}^{2+\a}}{\cA^{\a}(t)}+|\nabla \rho|_{\infty}^{2+\a}\right)+c_{\e}\cA^{\frac{\a}{\a-1}}(t),\nonumber\\
 |\nabla u|_{\infty}^3 &\leq \e \frac{|\nabla u|_{\infty}^{2+\a}}{\cA^{\a}(t)}+c_{\e}\cA^{\frac{3\a}{\a-1}}(t).
\end{align}
Therefore, we arrive at
\begin{align*}
\p_t \left[|\nabla \rho |_{\infty}^2+|\nabla u |_{\infty}^2\right] &\leq C |\nabla \rho|_{\infty} |\nabla e|_{\infty}+3 c_{\e}\cA^{\frac{\a}{\a-1}}(t)\\
& \quad  -(1-4 \e)\frac{|\nabla u|_{\infty}^{2+\a}}{\cA^{\a}(t)}-(1-3 \e)|\nabla \rho|_{\infty}^{2+\a} - (1-\e)D_{\a}(\nabla u)(x_{\star}) -D_{\a}(\nabla \rho)(x^{\star}). \nonumber
\end{align*}
To close the argument, we recall cf. \eqref{e:grade} that $|\nabla e (t)|_{\infty} \lesssim |\nabla \rho (t)|_{\infty}+ |\nabla u|_{L^1_{t} L^{\infty}_{x}}$ which forces to have an apriori control in time for the gradient of the velocity. We bypass this obstacle taking into account the fact that the ratio $q=e/\rho$ satisfies the transport equation \eqref{e:Dtq}. Consequently, we get
\begin{align}\label{e:rho+u+q}
\p_t \left[|\nabla \rho |_{\infty}^2+|\nabla u |_{\infty}^2 +|\nabla q|_{\infty}^2\right] &\lesssim C \left[|\nabla \rho|_{\infty} |\nabla e|_{\infty}+|\nabla u|_{\infty} |\nabla q|_{\infty}\right]+3 c_{\e}\cA^{\frac{\a}{\a-1}}(t)\\
& \quad  -(1-4 \e)\frac{|\nabla u|_{\infty}^{2+\a}}{\cA^{\a}(t)}-(1-3 \e)|\nabla \rho|_{\infty}^{2+\a} \nonumber \\ 
& \quad  - (1-\e)D_{\a}(\nabla u)(x_{\star}) -D_{\a}(\nabla \rho)(x^{\star}). \nonumber
\end{align}
In addition, by the definition of $q$ and the uniform bounds previously proved for density and entropy, we trivially have that
$|\nabla e|_{\infty}\lesssim |\nabla \rho|_{\infty}+|\nabla q|_{\infty}$. Therefore, the quadratic term of \eqref{e:rho+u+q} can be bounded directly as  
\[
|\nabla \rho|_{\infty} |\nabla e|_{\infty}+|\nabla u|_{\infty} |\nabla q|_{\infty}\lesssim |\nabla \rho|_{\infty}^2+|\nabla u|_{\infty}^2+|\nabla q|_{\infty}^2.
\]
Then,  we have obtained uniform bound  for $\rho, u, q \in L^{\infty}([0,T_0);\dot{W}^{1,\infty})$  by integration. This  fulfills the continuation criterion \eqref{e:BKMsing} and the proof of global existence of \eqref{e:CSHansatz} for $1<\a<2$ is complete.

\noindent {\sc Step 4.4: Strong Flocking.}
By the general result in multi-D proved in \cite{TT2014}, we have exponential alignment and flocking for any fat tail communication \eqref{e:ft}:
\begin{equation*}
\begin{split}
\mathcal{D}(t)&\leq \bar{\cD}<\infty \qquad \text{where }\quad \mathcal{D}(t):=\underset{\left\lbrace x, y \right\rbrace\in\Supp \rho(\cdot,t)}{\text{max}}|x- y | \\
\cA(t) &\leq \cA_0 e^{-\d t}.
\end{split}
\end{equation*}
Next, we complement this general result by a strong flocking statement of Theorem \ref{t:flock-singular}. This has so far been a 1D specific result, see \cite{ST2,ST3}, but we can extend it to the general multidimensional oriented flows and the use of the same entropy conservation. See our companion paper \cite{LS-uni1} for a similar resut in  the smooth kernel case.
The technical issue in applying the 1D strategy is that, again, $e$ only controls $\p_1 u$ and our goal is to extend such control to the full gradient and Hessian.\\

We already know from the previous  step that $|\nabla u|_{\infty}$ remains uniformly bounded. However, this argument does not provide a good quantitative estimate on $|\nabla u|_{\infty}$ to conclude flocking. We will seek more precise estimates with the fact that now we already know that  $|\nabla \rho|_{\infty}$ remains uniformly bounded. So, we go back to {\sc Step 4.2} to improve our estimate and finally get
\[
\p_t |\nabla u|_{\infty}^2 \lesssim C\, |\nabla u|_{\infty}^3  + c_{\e}\cA^{\a}(t)-(1-\e)\frac{|\nabla u|_{\infty}^{2+\a}}{\cA^{\a}(t)}.
\]
We absorb the cubic term simply by interpolation \eqref{e:cubic}, which again gets absorbed by cost of adding another $\cA^{\frac{3\a}{\a-1}}(t)$. In the end, we arrive at
\[
\p_t |\nabla u|_{\infty}^2 \lesssim -c \frac{|\nabla u|_{\infty}^{2+\a}}{\cA^{\a}(t)}+C\cA^{\beta}(t), \qquad \beta>0,
\]
and the result follows. In particular, it implies exponential rate of convergence to zero as $t\to\infty.$
Lastly, showing exponential decay of $|\nabla^2 u|_{\infty}$ follows similar estimates on the evolution of the norm $|\nabla^2 u|_{\infty}^{2}$, and will not be presented here for the sake of brevity. We refer to \cite{ST2} for full details in \nolinebreak 1D.

Now to establish strong flocking for the density we have that the velocity alignment goes to its natural limit $\bar{u}=\mathcal{P}/\mathcal{M}$. To do it, we pass to the moving frame $x-\bar{u}t$ and write the continuity equation in new coordinates. Then, $\tilde{\rho}(x,t):=\rho(x_{1}+t\bar{u},x_{2},\ldots,x_n,t)$ satisfies
$$\partial_{t}\tilde{\rho}(x,t)+\tilde{\rho}(x,t)\,\partial_{1} u(x_{1}+t\bar{u},x_{2},\ldots,x_n,t) +\partial_{1}\tilde{\rho}(x,t)(u(x_{1}+t\bar{u},x_{2},\ldots,x_n,t)-\bar{u})=0.$$
According to the established bounds we have that $|\partial_{t}\tilde{\rho}|_{\infty}=E(t)$, where in what follows $E(t)$ denotes a generic exponential decaying quantity. This shows that $\tilde{\rho}(\cdot,t)$ is Cauchy in $t$ in the metric of $L^{\infty}$. Hence, there exists a unique limiting state $\rho_{\infty}(\cdot)$ such that $|\tilde{\rho}(\cdot,t)-\rho_{\infty}(\cdot)|_{\infty}=E(t).$
Shifting $x_{1}$ this can be expressed in terms of $\rho(x,t)$ and $\bar{\rho}(x,t):=\rho_{\infty}(x_{1}-t\bar{u},x_{2},\ldots,x_n,t)$ as 
$|\rho(\cdot,t)-\bar{\rho}(\cdot)|_{\infty}= E(t).$
Since $\nabla \rho$ is uniformly bounded, this also shows that $\bar{\rho}$ is Lipschitz. Convergence
in $C^{\gamma}$ with $0<\gamma<1$ follows by interpolation.


\begin{thebibliography}{10}
	
	\bibitem{AC2019}
	Victor Arnaiz and \'Angel Castro.
	\newblock Singularity formation for the fractional euler-alignment system in
	1{D}.
	\newblock 2019.
	\newblock https://arxiv.org/abs/1911.08974.
	
	\bibitem{CV2012}
	Peter Constantin and Vlad Vicol.
	\newblock Nonlinear maximum principles for dissipative linear nonlocal
	operators and applications.
	\newblock {\em Geom. Funct. Anal.}, 22(5):1289--1321, 2012.
	
	\bibitem{CC2004}
	Antonio C\'{o}rdoba and Diego C\'{o}rdoba.
	\newblock A maximum principle applied to quasi-geostrophic equations.
	\newblock {\em Comm. Math. Phys.}, 249(3):511--528, 2004.
	
	\bibitem{CS2007a}
	Felipe Cucker and Steve Smale.
	\newblock Emergent behavior in flocks.
	\newblock {\em IEEE Trans. Automat. Control}, 52:852--862, 2007.
	
	\bibitem{CS2007b}
	Felipe Cucker and Steve Smale.
	\newblock On the mathematics of emergence.
	\newblock {\em Jpn. J. Math.}, 2(1):197--227, 2007.
	
	\bibitem{DMPW2019}
	Rapha\"{e}l Danchin, Piotr~B. Mucha, Jan Peszek, and Bartosz Wr\'{o}blewski.
	\newblock Regular solutions to the fractional {E}uler alignment system in the
	{B}esov spaces framework.
	\newblock {\em Math. Models Methods Appl. Sci.}, 29(1):89--119, 2019.
	
	\bibitem{DKRT2018}
	Tam Do, Alexander Kiselev, Lenya Ryzhik, and Changhui Tan.
	\newblock Global regularity for the fractional {E}uler alignment system.
	\newblock {\em Arch. Ration. Mech. Anal.}, 228(1):1--37, 2018.
	
	\bibitem{FK2019}
	Alessio Figalli and Moon-Jin Kang.
	\newblock A rigorous derivation from the kinetic {C}ucker-{S}male model to the
	pressureless {E}uler system with nonlocal alignment.
	\newblock {\em Anal. PDE}, 12(3):843--866, 2019.
	
	\bibitem{HL2009}
	Seung-Yeal Ha and Jian-Guo Liu.
	\newblock A simple proof of the {C}ucker-{S}male flocking dynamics and
	mean-field limit.
	\newblock {\em Commun. Math. Sci.}, 7(2):297--325, 2009.
	
	\bibitem{HT2008}
	Seung-Yeal Ha and Eitan Tadmor.
	\newblock From particle to kinetic and hydrodynamic descriptions of flocking.
	\newblock {\em Kinet. Relat. Models}, 1(3):415--435, 2008.
	
	\bibitem{KNV2007}
	A.~Kiselev, F.~Nazarov, and A.~Volberg.
	\newblock Global well-posedness for the critical 2{D} dissipative
	quasi-geostrophic equation.
	\newblock {\em Invent. Math.}, 167(3):445--453, 2007.
	
	\bibitem{LS-uni1}
	Daniel Lear and Roman Shvydkoy.
	\newblock Existence and stability of unidirectional flocks in hydrodynamic
	{E}uler {A}lignment systems.
	\newblock 2019.
	\newblock https://arxiv.org/abs/1911.10661.
	
	\bibitem{LS-entropy}
	Trevor~M. Leslie and Roman Shvydkoy.
	\newblock On the structure of limiting flocks in hydrodynamic {E}uler
	{A}lignment models.
	\newblock {\em Math. Models Methods Appl. Sci.}, 29(13):2419--2431, 2019.
	
	\bibitem{SS2016}
	Russell~W. Schwab and Luis Silvestre.
	\newblock Regularity for parabolic integro-differential equations with very
	irregular kernels.
	\newblock {\em Anal. PDE}, 9(3):727--772, 2016.
	
	\bibitem{Shv2018}
	Roman Shvydkoy.
	\newblock Global existence and stability of nearly aligned flocks.
	\newblock {\em J. Dynam. Differential Equations}, 31(4):2165--2175, 2019.
	
	\bibitem{ST1}
	Roman Shvydkoy and Eitan Tadmor.
	\newblock Eulerian dynamics with a commutator forcing.
	\newblock {\em Transactions of Mathematics and Its Applications}, 1(1):tnx001,
	2017.
	
	\bibitem{ST2}
	Roman Shvydkoy and Eitan Tadmor.
	\newblock Eulerian dynamics with a commutator forcing {II}: {F}locking.
	\newblock {\em Discrete Contin. Dyn. Syst.}, 37(11):5503--5520, 2017.
	
	\bibitem{ST3}
	Roman Shvydkoy and Eitan Tadmor.
	\newblock Eulerian dynamics with a commutator forcing {III}. {F}ractional
	diffusion of order {$0<\alpha<1$}.
	\newblock {\em Phys. D}, 376/377:131--137, 2018.
	
	\bibitem{S2012}
	Luis Silvestre.
	\newblock H\"older estimates for advection fractional-diffusion equations.
	\newblock {\em Ann. Sc. Norm. Super. Pisa Cl. Sci. (5)}, 11(4):843--855, 2012.
	
	\bibitem{TT2014}
	Eitan Tadmor and Changhui Tan.
	\newblock Critical thresholds in flocking hydrodynamics with non-local
	alignment.
	\newblock {\em Philos. Trans. R. Soc. Lond. Ser. A Math. Phys. Eng. Sci.},
	372(2028):20130401, 22, 2014.
	
	\bibitem{Tan2017}
	Changhui Tan.
	\newblock Singulaity formation for a fluid mechanics model with nonlocal
	velocity.
	\newblock 2017.
	\newblock https://arxiv.org/abs/1708.09360.
	
\end{thebibliography}

\end{document}